\documentstyle{article}

\font\Bbb=msbm10 scaled 1000 

\newcommand{\pl}{\partial_}
\newcommand{\ve}{\varepsilon_}
\newcommand{\ol}{\overline}
\newcommand{\bqq}{\begin{equation} \label}
\newcommand{\eeq}{\end{equation}}
\newcommand{\rr}{{\mbox{\Bbb R}}}
\newcommand{\cc}{{\mbox{\Bbb C}}}

\newcommand{\tmc}{$TM \otimes_{\mbox{\Bbb R}}\mbox{\Bbb C}$}
\newcommand{\dmin}{(\partial_1- \partial_{\overline{1}})}

\newcommand{\appsection}{\addtocounter{section}{1}
           \setcounter{equation}{0}\section*{Appendix \Alph{section}}}

\newtheorem{The}{Theorem}
\newtheorem{Lem}{Lemma}
\newtheorem{Nte}{Note}

\begin{document}
\hsize=15cm
\textheight=24cm
\addtolength{\topmargin}{-3cm}
\mathsurround=2pt

\begin{center}
{\LARGE $H$-projectively-equivalent
\medskip

four-dimensional K\"ahler manifolds}
\bigskip \bigskip

{\Large Dmitry A. Kalinin}

Department of General Relativity \\
and Gravitation Kazan State University \\
18 Kremlyovskaya Ul. Kazan 420008 Russia\\
E-mail: {\tt Dmitry.Kalinin{\it @}ksu.ru}
\end{center}

\bigskip \bigskip \bigskip

\begin{abstract}
The paper is devoted to the investigation of four-dimensional
K\"ahler manifolds admitting non-affine $H$-projective mappings. We
find all such manifolds which are non-Einstein. In the paper also K\"ahler
manifolds admitting infinitesimal $H$-projective transformations are
de\-ter\-mi\-ned. It is proved that the class of K\"ahler manifolds
admitting $H$-projective mappings includes generalized
equidistant K\"ahler manifolds. The explicit expression for the metrics of
Ricci-flat four-dimensional generalized equ\-i\-dis\-t\-ant
manifolds is given.
\end{abstract}
\bigskip

\section{Introduction}
K\"ahler manifolds was introduced by P.~A.~Shirokov~\cite{shir1} and
E.~K\"ahler \cite{kahl} in the first part of our century. Since that
time they gained applications in a wide variety of fields both in
mathematics and physics. These fields are algebraic geometry,
dynamical systems, general relativity theory, superstring theory and
nonlinear $\sigma$-models, quantum theory {\it etc.}
\cite{al-free,amka2,boraj,flah,kalinin,kalinin1}.

In 1954 Otsuki and Tashiro \cite{otas} have generalized the projective
mappings of (pseudo) Riemannian manifolds \cite{am1,sin1} by
introducing the notion of $H$-projective mappings of K\"ahler
manifolds. It is well-known that investigation of dynamical systems of
second order with common trajectories can be reduced to solving the
problem of finding pseudo Riemannian (and affine connected) manifolds
admitting projective mapping. Similarly, the investigation of
$H$-projective mappings of K\"ahler manifolds is closely related to
the problem of trajectory equivalency of dynamical systems of second
order with non-zero external forces of special
form \cite{kalinin,kalinin1}.

To the present time wide variety K\"ahler manifolds {\it not
admitting} $H$-projective mappings is found. These are equidistant
\cite{mik-ekv}, generalized symmetric and recurrent manifolds
different from manifolds of constant holomorphic sectional
curvature \cite{sin1}. At the same time, some general methods of
finding $H$-projective mappings for given K\"ahler manifold were also
developed \cite{sin1,skm}.

However, many essential problems to be solved remain in the theory of
$H$-projective mappings. Above all should be referred the problem of
finding Riemannian metrics and connections for K\"ahler manifolds
admitting non-affine $H$-projective mappings. To the present this
problem is unsolved ever in the case of lower dimensions. Some
approaches to its solution was proposed earlier \cite{amka1} by the
author.

The present paper is devoted to finding of four-dimensional K\"ahler
manifolds admitting non-affine $H$-projective mappings. We solve this
problem for non-Einstein manifolds as well as in the case of
infinitesimal transformations. We proved that non-Einstein K\"ahler
manifolds admitting $H$-projective mappings are generalized
equidistant manifolds. Moreover, it is possible to prove that
generalized equidistant K\"ahler manifolds admit $H$-progective
mappings in general case. In the last part of the paper Einstein and
Ricci-flat generalized equidistant K\"ahler manifolds having special
interest from the physical point of wiev are considered. explicit
expression for the metrics of Ricci-flat four-dimensional generalized
equidistant manifolds is found. This paper is corrected and extended
version of \cite{izv2} (see also \cite{diss}).


\section{Differential geometry of K\"ahler manifolds}


Let us start from reminding some relevant facts on differential
geometry of K\"ahler manifolds \cite{kobn,sin1}.

An $2n$-dimensional ($n>1$) smooth manifold $M$ is called to be {\it
almost complex} if the {\it almost complex structure} $J_p:T_M \to
T_M$, $J^2=-{\rm id}|_{TM}$ is defined in its tangent bundle. An
(1,2)-ten\-sor field $N$ on $M$ defined by the formula
     $$
N(X,Y) = 2 ([JX,JY]- [X,Y] - J[X,JY] - J[JX,Y]),
     $$
for any vector fields $X,Y$ is called {\it torsion} of almost complex
structure $J$. If $N$ vanishes then $J$ is called {\it complex
structure}. In this case $(M,J)$ is called to be {\it integrable
manifold}.

Let $M$ be integrable almost complex manifold with complex structure
$J$. According to Newlander-Nirenberg theorem \cite{kobn}, there
exists the unique complex manifold $M^c$ coinciding with $M$ as
topological space and such that its complex analytic structure induces
the complex structure $J$ and the structure of differential manifold
on $M$.

The tangent bundle $TM^c$ is {\Bbb C}-linear isomorphic to the bundle $TM$
with the structure of complex bundle induced by the operators $J_p$,
$p\in M$ so that there is a canonical {\Bbb C}-linear bundle
isomorphism
  \bqq{1a}
TM \otimes_\rr\cc \cong TM^c  \oplus \ol{TM^c}
  \eeq
where $TM\otimes_\rr\cc$ is complexification of tangent bundle $TM$
and bar denotes the complex conjugation.

Let $(U, z^{\alpha})$, $\alpha=1,...,n$ be a chart on $M^c$.
If $M$ is integrable almost complex manifold
corresponding to $M^c$ then we shall say that $(U,z^\alpha,
\ol{z^\alpha})$, $\alpha=1,...,n$ is {\it complex chart} on $M$.
Because of isomorphism (\ref{1a}) the vector fields $\partial_\alpha
\equiv \partial / \partial z^{\alpha}, \partial_{\ol \alpha} \equiv
\partial / \partial \overline{z^{\alpha}}$, $\alpha=1,...,n$ define a
basis in \tmc. Any real tensor field $T$ on $M$ can be uniquely
extended to the smooth field of elements of tensor algebras
     $$
\tilde{\bf T}_p M \equiv \bigoplus \limits^{\infty}_{k_i=1} ( (T_p
M^c)^{\otimes k_1}\otimes (\ol{T_p M^c})^{\otimes k_2} \otimes (T^*_p
M^c)^{\otimes k_3}\otimes (\ol{T^*_p M^c})^{\otimes k_4}).
     $$
As well as tensor $T$, this extension also is called tensor field and
is denoted by the same letter $T$. The expansion of $T$ with respect
to the basis $(\pl{\alpha}$, $\pl{\ol\alpha})_{\alpha=1,...,n}$ has
the form
     $$
T = T^{i_1 ... i_{r}}_{j_1 ... j_s} \partial_{i_1} \otimes ... \otimes
\partial_{i_{r}} \otimes dz^{j_1} \otimes ... \otimes dz^{j_s}, \quad
T^{\ol{i_1}...\ol{i_{r}}}_{\ol{j_1}... \ol{j_s}} = \ol{T^{i_1 ...
i_{r}}_{j_1 ... j_s}}.
     $$
Here Latin indices $i,j,k,... = 1,...,2n$ run over the sets of bared
$(\ol \alpha, \ol \beta, \ol \gamma,...)$ and unbarred $(\alpha,
\beta,\gamma,...)$ Greek indices which are varied from $1$ to $n$. The
functions $T^{i_1 ... i_r}_{j_1 ... j_s}$ are called to be {\it
components of tensor field $T$ in coordinates} $(z^\alpha,
\ol{z^\alpha})$.

In particular, complex structure $J$ can be uniquely extended to
{\Bbb C}-linear endomorphism in \tmc. The action of complex structure
on the elements of coordinate basis $(\partial_{\alpha}$,
$\partial_{\bar \alpha})$, $\alpha=1,...,n$ is defined by the
equalities $J\partial_\alpha= i\partial_\alpha$, $J\partial_{\ol
\alpha}= -i\partial_{\ol \alpha}$.

Let us call {\it holomorphic transformation} the coordinate
transformation of the form $z'^{\alpha}=w^{\alpha}(z)$,
${z'}^{\ol\alpha}=\ol{w^\alpha (z)}$ where $w^{\alpha}(z)$ are complex
analytic functions. Let $X = \xi^i\partial_i$, $\xi^{\ol\mu} =
\ol{\xi^\mu}$ be a real vector field on $M$. If the Lie derivative
$L_X J$ is equal to zero then $X$ is called to be {\it holomorphic
vector field}. This condition in a complex chart $(U,
z^{\alpha},z^{\ol\alpha})_{\alpha=1,...,n}$ reduces to the
equalities $\partial_{\ol\nu} \xi^\mu= \partial_{\nu} \xi^{\ol\mu}
=0$, $\mu,\nu =1,...,n$. Using holomorphic coordinate transformations
any holomorphic vector field can be reduced to the form
$X=\pl 1 + \pl{\ol 1}$ \cite{skm}.

An integrable almost complex manifold $M$ is called to be {\it
K\"ahler manifold} if a pseudo Riemannian metric $g$ is defined on $M$
satisfying the following conditions \cite{kobn,sin1}
     \bqq{sogl1}
g(JX,JY)=g(X,Y), \qquad \nabla J = 0,\qquad (X,Y \in TM).
     \eeq
Here $\nabla$ is Levi-Civita connection of the metric $g$. The
bilinear form
    \bqq{333}
\Omega(X,Y)=g(JX,Y)
    \eeq
is called {\it fundamental} 2-{\it form} of K\"ahler manifold $M$.
>From Eqs.~(\ref{sogl1}), (\ref{333}) and condition $J^2 = -{\rm
id}|_{TM}$ it follows that $\Omega$ is closed 2-form: $d\Omega=0$.

Let $(U,z,\ol z)$ be complex chart on $(M,g)$ and $\pl{\alpha}$,
$\pl{\ol\alpha}$, $\alpha=1,...,n$ are coordinate vector fields.
Then the components of the metric $g$, complex structure $J$ and
fundamental
2-form $\Omega$ in this chart are defined by the folowing conditions
   \bqq{mherm}
g_{\alpha\ol\beta}=\ol {g_{\ol\alpha\beta}},\qquad
g_{\alpha\beta}=g_{\ol\alpha\ol\beta}=0,\qquad
g^{\alpha\ol\beta}=\ol{g^{\ol\alpha\beta}}, \qquad
g^{\alpha\beta}=g^{\ol\alpha\ol\beta}=0,
   \eeq
   \bqq{comJ}
J^\alpha_{\beta}=- J_{\ol\beta}^{\ol\alpha}=i\delta^\alpha_\beta,\qquad
J^\alpha_{\ol\beta} = J_{\beta}^{\ol\alpha} = 0,
   \eeq
   \bqq{com-omega}
\Omega_{\alpha\ol\beta}=-i g_{\alpha\ol\beta}, \qquad
\Omega_{\ol\alpha\beta}=- \ol{\Omega_{\alpha\ol\beta}},\qquad
\Omega_{\alpha\beta}=\Omega_{\ol\alpha\ol\beta}.
\eeq
In the same coordinates the condition $d\Omega=0$ takes the form
   \bqq{zamkn}
\partial_\alpha g_{\beta\ol\gamma}=\partial_\beta g_{\alpha\ol\gamma},
\qquad
\pl{\ol\alpha} g_{\ol\beta\gamma} = \pl{\ol\beta}g_{\ol\alpha\gamma}.
   \eeq
>From here it follows that in $U$ exists a real-valued function $\Phi$
obeying
   \bqq{kahl}
g_{\alpha\ol\beta}=\partial_{\alpha}\partial_{\ol\beta}\Phi.
   \eeq
This function is called {\it K\"ahler potential} of the metric $g$ and
is defined up to the transformations
  \bqq{gau}
\Phi' = \Phi +f(z) +\ol{f(z)}
  \eeq
where $f(z)$ is an appropriate holomorphic function. Such
transformations are called {\it gauge transformations} of K\"ahler
potential \cite{izv2}.

>From Eqs.~(\ref{mherm})--(\ref{kahl}) it follows that the only
non-vanishing
Christoffel symbols are
\bqq{cris}
\Gamma^\alpha_{\beta\nu}= \ol{\Gamma^{\ol\alpha}_{\ol\beta
\ol\nu}}= g^{\alpha\ol\mu}\partial_\beta g_{\ol\mu\nu}
\eeq
and non-zero components of Riemann and Ricci curvature tensors are
defined by the conditions
    \bqq{riem}
R^{\alpha}_{\beta\mu\ol\nu} = \ol{ R^{\ol\alpha}_{\ol\beta\ol\mu\nu} }
= -R^{\alpha}_{\beta\ol\nu\mu} =
-\ol{R^{\ol\alpha}_{\ol\beta\nu\ol\mu} } = -\partial_{\ol\nu}
\Gamma^{\alpha}_{\beta\mu},
    \eeq
    \bqq{ric}
R_{\alpha\ol\beta} = \partial_{\alpha}\partial_{\ol\beta} \: \mbox{ln}
(\det (g_{\mu\ol\nu})), \qquad R_{\alpha\ol\beta} =
\ol{R_{\ol\alpha\beta}}.
    \eeq


\section{H-projective mappings of K\"ahler manifolds}


A smooth curve $\gamma: [0,1] \to M :t \mapsto x_t$ on K\"ahler
manifold $M$ of real dimension $2n>2$ is called to be {\it $H$-planar
curve} if its tangent vector $\chi \equiv dx/dt$ obeys the equations
     $$
\nabla_{\chi} \chi = a(t)\chi +  b(t) J(\chi)
     $$
where $a(t)$ and $b(t)$ are functions of the parameter $t$.

Let us consider two K\"ahler manifolds $M$, $M'$ with metrics $g$,
$g'$ and complex structures $J$, $J'$. A diffeomorphism $f: M \to M'$
is called $H$-{\it projective mapping} if for any $H$-planar curve
$\gamma$ in $M$ the curve $f\circ \gamma$ is $H$-planar curve in $M'$.
Any $H$-projective mapping preserves the complex structure $J$
\cite{skm}, i.e. $f_* \circ J =J' \circ f_*$ where $f_*$ is
the differential of $f$. From here it follows that it is possible to
choose such complex charts in $M$ and $M'$ that the corresponding
points $p \in M$ and $p' = f(p) \in M'$ have the same coordinates and
components of complex structures $J$, $J'$ coincide and have the
form (\ref{comJ}). Such coordinate systems in $M$, $M'$ is said to be
{\it corresponding complex coordinates with respect to the mapping}
$f$.

Necessarily and sufficient condition of $f$ to be $H$-projective
mapping can be expressed in corresponding coordinates by the equation
\cite{sin1}
   \bqq{sh-bis1}
\Gamma'^{i}_{jk} - \Gamma^{i}_{jk} = 2 \delta ^{i}_{(j} \psi _{,k)} -
2 \psi _{,l} J ^{l}_{(j} J^{i}_{k)}
   \eeq
where $\psi$ is a real valued function, $\Gamma'^{i}_{jk}$,
$\Gamma^{i}_{jk}$ are Christoffel symbols of the metrics $g$, $g'$ and
$\psi_i = \psi_{,i} \equiv \pl i \psi$ where comma denotes the
covariant derivative with respect to $g$. If, in particular,
$\psi_k=0$ and $\psi=$const, then $H$-projective
mapping preserves Riemannian connection and is called {\it affine}
mapping. The Eq.~(\ref{sh-bis1}) is equivalent to the following
equation \cite{amka1}
         $$
\nabla g'(Y,Z,W) = 2g'(Y,Z) W\psi + g'(Z,W) Y\psi +g'(Y,W)Z\psi
         $$
         \bqq{hp1}
- (JY)(\psi) g'(Z,JW) - (JZ)(\psi) g'(Y,JW),\qquad (Y,Z,W \in TM).
         \eeq
Setting here $Y=\partial_\alpha$, $Z=\partial_{\ol\beta}$ and $W =
\partial_\gamma$, with the help of (\ref{mherm}), (\ref{comJ}) we find
         $$
g'_{\alpha\ol\beta,\gamma} = 2 g'_{\alpha\ol\beta}\psi_\gamma + 2
g'_{\gamma\ol\beta}\psi_\alpha, \qquad g'_{\alpha\beta,\gamma} =
g'_{\alpha\beta,\ol\gamma}=0,
         $$
         \bqq{hp2}
g'_{\alpha\ol\beta,\ol\gamma} = 2 g'_{\alpha\ol\beta}\psi_{\ol\gamma}
+ 2 g'_{\alpha\ol\gamma}\psi_{\ol\beta}, \qquad
g'_{\ol\alpha\ol\beta,\gamma} = g'_{\ol\alpha\ol\beta,\ol\gamma}=0.
         \eeq
Using {\it Sinyukov's $\Gamma$-transformation} \cite{sin1}
         \bqq{sin}
a_{\alpha\ol\beta}= \ol{a_{\ol \alpha\beta}} =
e^{2\psi}g'^{\lambda\ol\mu}g_{\alpha\ol\mu} g_{\lambda\ol\beta},
\qquad a_{\alpha\beta}=a_{\ol\alpha\ol\beta}=0, \qquad
g^{\alpha\ol\beta}= e^{-2\psi } a^{\alpha\ol\beta}
         \eeq
where
         $$
a^{\alpha\ol\beta} = a_{\mu\ol\lambda}g^{\alpha\ol\lambda}
g^{\mu\ol\beta}, \qquad (g'^{\alpha\ol\beta})= (g'_{\alpha\ol\beta})^{-1},
         $$
we can write Eq.~(\ref{hp2}) in the form
         \bqq{hpa}
a_{\alpha\ol\beta,\gamma} = \lambda_\alpha g_{\gamma\ol\beta}, \qquad
a_{\alpha\ol\beta,\ol\gamma} = \lambda_{\ol\beta} g_{\gamma\ol\alpha}
         \eeq
where
         $$
\lambda_\alpha = \ol{\lambda_{\ol\alpha}} =
- 2 \psi_\nu e^{2 \psi}g'^{\nu\ol\mu} g_{\alpha\ol\mu}.
         $$
Contracting (\ref{hpa}) with $g^{\alpha\ol\beta}$, we obtain
         \bqq{lambda}
\lambda_\gamma = \frac{1}{2} \pl\gamma (g^{ij} a_{ij}) =
\partial_{\gamma} (a_{\alpha\ol\beta} g^{\alpha\ol\beta}),
         \qquad
\lambda_{\ol\gamma} = \frac{1}{2}\pl{\ol\gamma} (g^{ij} a_{ij}) =
\partial_{\ol\gamma} (a_{\alpha\ol\beta} g^{\alpha\ol\beta}) =
\lambda_{\ol\gamma}.
         \eeq
>From here it follows, that $\lambda_i$ is a gradient field, i.e. there
exists a real function $\lambda$ such that $\lambda_i dz^i=d\lambda$
and $\lambda_i = \lambda_{,i}$.

The integrability conditions of Eq.~(\ref{hpa}) follows from the Ricci
identity
         \bqq{ui-ij}
2a_{kl,[ij]} =a_{sl}R^s_{kij} + a_{ks} R^s_{lij}.
         \eeq
For $(ijkl)= (\gamma\ol\nu\alpha\ol\beta),(\gamma\nu\alpha\ol\beta)$
using (\ref{riem}) we find the following conditions
         \bqq{ui-1}
a_{\mu\ol\beta} R^{\mu}_{\alpha\gamma\ol\nu}+
a_{\alpha\ol\mu}R^{\ol\mu}_{\ol\beta\gamma\ol\nu}=
g_{\gamma\ol\beta}\lambda_{\alpha,\ol\nu}-
g_{\alpha\ol\nu}\lambda_{\ol\beta,\gamma},
         \eeq
         \bqq{ui-2}
g_{\gamma\ol\beta}\lambda_{\alpha,\nu}-
g_{\nu\ol\beta}\lambda_{\alpha,\gamma}=0.
         \eeq
The remaining integrability conditions hold identitically or can be
obtained from (\ref{ui-1}) and (\ref{ui-2}) by complex conjugation.
Contracting (\ref{ui-1}) with $g^{\alpha\ol\nu}$, with the use of the
formula $R^{\alpha}_{\beta\gamma\ol\nu}g^{\beta\ol\nu} =
-R^{\alpha}_{\gamma}$, we receive
         $$
- a_{\mu\ol\beta} R^{\mu}_{\gamma}+
a_{\alpha\ol\mu}{R^{\ol\mu}_{\ol\beta\gamma}}^{\alpha}=
g_{\gamma\ol\beta} g^{\alpha\ol\nu} \lambda_{\alpha,\ol\nu}- n
\lambda_{\ol\beta,\gamma}.
         $$
>From here using (\ref{lambda}) and $a_{\alpha\ol\mu}
{R^{\ol\mu}_{\ol\beta\gamma}}^{\alpha}= a_{\alpha\ol\mu}
{R^{\alpha}_{\gamma\ol\beta}}^{\ol\mu}$, it is easy to derive that
$a_{\mu\ol\beta} R^{\mu}_{\gamma}- a_{\gamma\ol\mu}
R^{\ol\mu}_{\ol\beta}=0$. Contracting this equation with
$g^{\ol\beta\nu}$, we find
         \bqq{zvez}
a_{\mu}^{\nu} R^{\mu}_{\gamma}- a_{\gamma}^{\mu} R^{\nu}_{\mu}=0.
         \eeq
Note, that as a consequence of (\ref{mherm}), (\ref{hpa}) we have
$a^{\ol\alpha}_{\beta} = a^\alpha_{\ol\beta}=0$. Contracting
(\ref{ui-2}) with $g^{\gamma\ol\beta}$ one can find
$(n-1)\lambda_{\alpha,\nu} =0$ which means that $\lambda_{\alpha,\nu}
=0$ and $\lambda^{\alpha}_{,\ol\nu} =0$, or, because of
$\Gamma^\alpha_{\ol\nu i} =0$,
   \bqq{hol-lam}
\partial_{\ol\nu} \lambda^{\alpha} =0, \qquad \partial_{\nu}
\lambda^{\ol\alpha} =0.
    \eeq
So, we come to conclusion that $\Lambda = \lambda^i \partial_i$ is
holomorphic vector field. Using holomorphic coordinate transformations
this vector field can be reduced to the form
    \bqq{hol-form}
\Lambda = \pl1 +\pl{\ol 1}, \qquad \lambda^\alpha = \delta^\alpha_1,
\qquad \lambda^{\ol\alpha} = \delta^\alpha_1.
    \eeq

    \begin{The}  \label{sh-th1}
Let $f$ be non-affine $H$-projective mapping of K\"ahler manifold
$(M,g)$ on K\"ahler manifold $(M',g')$. Let also $d\lambda =
\lambda_{\alpha} dz^\alpha+ \lambda_{\ol\alpha} dz^{\ol\alpha}$ be the
exact 1-form defined by Eqs.~{\rm (\ref{hp2}) -- (\ref{lambda})}. Then
the real vector field $J\Lambda = i\lambda^{\alpha} \pl\alpha -
i\lambda^{\ol\alpha} \pl{\ol\alpha}$ is infinitesimal isometry of
$M$, i.e. the following Killing equations hold: $L_{J\Lambda} g =0$.
     \end{The}

\noindent {\bf Proof:} Using Eqs.~(\ref{cris}), (\ref{lambda}) and
(\ref{hol-lam})
     $$
-i\lambda_{\ol\beta,\alpha} +i\lambda_{\alpha,\ol\beta} = -i
\pl{\alpha}\pl{\ol\beta} (a_{\nu\ol\mu} g^{\nu\ol\mu})+
 i \pl{\ol\beta}\pl{\alpha} (a_{\nu\ol\mu} g^{\nu\ol\mu})=0,
     $$
     $$
i\lambda_{\beta,\alpha} + i\lambda_{\alpha,\beta} = 0, \qquad
-i\lambda_{\ol\beta,\ol\alpha} - i\lambda_{\ol\alpha,\ol\beta} = 0,
     $$
or
     \bqq{killing}
L_\mu g_{ij} \equiv \mu_{i,j} + \mu_{j,i} =0, \qquad \mu=J\Lambda
     \eeq
where $\mu_i =g_{il} \mu^l$. Since $\Lambda$ and $\mu$ don't vanish
(since $f$ is non-affine mapping) and $\mu^i\pl i = J\Lambda$,
then $J \Lambda$ is infinitesimal isometry. {\it Q.E.D.}

If we reduce $\Lambda$ to (\ref{hol-form}), then the Killing
equations take the form
    \bqq{nondep1}
(\pl1 -\pl{\ol 1}) g_{\alpha\ol\beta} = 0.
    \eeq

    \begin{Lem} \label{sh-lm1}
If a K\"ahler manifold $(M,g)$ admits infinitesimal isometry
$J\Lambda$ where $\Lambda$ is defined by {\rm (\ref{hol-form})},
then K\"ahler potential of $g$ can be reduced to the following form
      \bqq{sh-bis2}
\Phi = \Phi (z^1 +z^{\ol 1}, z^2, z^{\ol 2}, ...).
      \eeq
      \end{Lem}

\noindent{\bf Proof:} Using (\ref{kahl}), we find from (\ref{nondep1})
      \bqq{nondep2}
(\pl 1  -\pl{\ol  1})  g_{\alpha\ol\beta}  =  \pl{\alpha}\pl{\ol\beta}
\dmin\Phi =0.
      \eeq
Hence, $\dmin \Phi= f(z) + h(\ol z)$ where $f$ is a holomorphic
function and $h$ is an antiholomorphic function. Similarly, because of
reality of the K\"ahler potential we have $(\pl 1 - \pl{\ol 1}) \Phi = -
\ol{(\pl 1 -\pl{\ol 1}) \Phi}$ and $h(\ol z) =- \ol{f(z)}$. Let us
change the K\"ahler potential using gauge transformations (\ref{gau})
      $$
\Phi = \Phi' +\int f(z) dz^1 + \ol{\int f(z) d z^1}= \int f(z)dz^1
-\int h(\ol z) d \ol{z^1}.
      $$
Substituting this expression in (\ref{nondep2}), we obtain $(\pl 1
-\pl{\ol 1}) \Phi' =0$. From here, omitting the primes we find $\Phi= \Phi
(z^1 + z^{\ol 1}, z^2, z^{\ol 2}, ...)$. {\it Q.~E.~D.}

In the following, if the conditions of Lemma~\ref{sh-lm1} are
satisfied, then we will suppose that K\"ahler potential is reduced to
the form (\ref{sh-bis2}).

Let a K\"ahler manifold $(M,g)$ admits non-affine $H$-projective
mapping, then, according to Theorem~\ref{sh-th1} and
Lemma~\ref{sh-lm1} it admits infinitesimal isometry $J\Lambda =i\dmin$
and its K\"ahler potential can be reduced to the form (\ref{sh-bis2}).
At the same time (\ref{hpa}) yield
  \bqq{hpa-mod1}
a^{\alpha}_{\beta,\ol\gamma}= \pl{\ol\gamma} a^{\alpha}_{\beta} =
\delta^{\alpha}_{1} g_{\beta\ol\gamma}, \qquad a^{\alpha}_{\beta} =
g^{\alpha\ol\sigma} a_{\beta\ol\sigma},
  \eeq
  \bqq{hpa-mod2}
a^{\alpha}_{\beta,\gamma}= \lambda_{,\beta} \delta^\alpha_\gamma
  \eeq
and complex conjugated equations. From (\ref{hpa-mod2}) it follows
$\lambda = a^\alpha_\alpha = a_{\ol\alpha}^{\ol\alpha} =
a_{\alpha\ol\beta} g^{\alpha\ol\beta}$. Integrating
Eq.~(\ref{hpa-mod1}), we find
  \bqq{hpa-sol1}
a^\alpha_\beta = \ol{a^{\ol\alpha}_{\ol\beta}} = \delta^\alpha_1
\pl\beta \Phi + h^\alpha_\beta
      \eeq
where $h^\alpha_\beta$ are holomorphic functions. Hence, $\lambda =
\pl 1 \Phi + h^\alpha_\alpha$. Because of reality of $\lambda$ and
$\pl{1} \Phi$
    \bqq{rho1}
h^\alpha_\alpha =h^{\ol\alpha}_{\ol\alpha} \equiv n\rho =
\mbox{const}, \qquad \lambda = \pl 1 \Phi +n\rho.
     \eeq

\begin{Nte}\label{note1}
Here we have used the fact that a holomorphic function is real iff it
is constant.
\end{Nte}
Substituting (\ref{rho1}) in (\ref{hpa-mod2}), we have
    \bqq{hpa-mod3}
a^{\alpha}_{\beta,\gamma}= g_{\beta \ol 1} \delta^\alpha_\gamma.
    \eeq


\section{Four-dimensional non-Einstein manifolds}


Let $(M_{4},g)$ be a non-Einstein ($R^i_j\neq\kappa \delta^i_j$)
K\"ahler manifold of dimension ${\mbox{dim}}_\rr M_{4} =4$
with K\"ahler potential $\Phi$. Let $M$ admits non-affine
$H$-projective mapping on a K\"ahler manifold $(M'_{4},g')$ and let
$a$ be the tensor field be defined by Eq.~(\ref{sin}). We can introduce
tensor field $b^i_j= L_{J\Lambda}a^i_j$ where $J \Lambda$ is corresponding
infinitesimal isometry (see Theorem~\ref{sh-th1}).

If we choose local complex coordinates such that
   \bqq{lam-jlam}
\Lambda=\lambda^i\pl i=\pl 1+\pl{\ol 1}, \qquad J \Lambda = i\dmin.
   \eeq
then, according to (\ref{hpa-sol1}) and Lemma~\ref{sh-lm1},
   \bqq{hpa-sol2}
a^\alpha_\beta = \ol{a^{\ol\alpha}_{\ol\beta}} = \delta^\alpha_1
\pl\beta \Phi + f^\alpha_\beta (z^1,z^2) +\rho \delta^\alpha_\beta,
\qquad f^\alpha_\alpha= f^{\ol\alpha}_{\ol\alpha} = 0,
   \eeq
   \bqq{potl}
\Phi = \Phi(z^1 +z^{\ol 1},z^2,z^{\ol 2})
   \eeq
where $f^\alpha_\beta \equiv h^\alpha_\beta-\rho \delta^\alpha_\beta$
are holomorphic functions. From here
    \bqq{hpab-sol}
b^\alpha_\beta=\ol{b^{\ol\alpha}_{\ol\beta}}=i\dmin a^\alpha_\beta =
i\pl 1 f^\alpha_\beta, \quad b^{\ol\alpha}_\beta=
b^\alpha_{\ol\beta}=0, \quad
b^\sigma_\sigma=b^{\ol\sigma}_{\ol\sigma}=0,
    \eeq

Admissible coordinate and gauge transformations which don't change
the form of vector field $\Lambda=\pl 1 +\pl{\ol 1}$ and the form
(\ref{potl}) of K\"ahler potential are
      \bqq{free1}
{z'}^1 = z^1+l(z^2), \qquad {z'}^2 = m(z^2),
      \eeq
      \bqq{free-gau}
\Phi'= \Phi +r\cdot (z^1 + \ol{z^1}) + p(z^2) + \ol{p(z^2)}, \qquad
r\in \rr
      \eeq
where $l$, $m$, $p$ are holomorphic functions depending on $z^2$.
Taking the Lie derivative $J\Lambda$ from both parts of (\ref{zvez}),
we find
      \bqq{zvez-b}
b_{\mu}^{\nu} R^{\mu}_{\gamma}- b_{\gamma}^{\mu} R^{\nu}_{\mu}=0,
     \eeq
or
     $$
b_2^1 R^2_1 - b^2_1 R_2^1 = 0,
     $$
     \bqq{zvez-bb}
2 b_1^1 R^1_2 + b^1_2 (R^2_2 - R^1_1) = 0,
     \eeq
     $$
2 b^1_1 R^2_1 + b^2_1 (R_2^2 - R^1_1) = 0.
     $$

Using this formulas it is possible to prove the following

\begin{Lem}\label{sh-lm2}
If a non-Einstein four-dimensional K\"ahler manifold $M_4$ admits
non-affine $H$-projective mapping, then in a neighborhood of any
point $p \in M_4$ exist complex coordinates in which the following
relations hold
   \bqq{sh48}
a^\alpha_\beta =\delta_1^\alpha \pl\beta\Phi +f_\beta^\alpha(z^2)
+\rho \delta_\beta^\alpha, \qquad \dmin\Phi =0.
   \eeq
   \end{Lem}

We have placed the proof, which is rather long and technical, in
Appendix~A so as not interrupt exposition.

Admissible coordinate and gauge transformations not changing
Eq.~(\ref{sh48}) are defined by the formulas (\ref{free1}) and
(\ref{free-gau}). Using these transformations one can reduce
$f^\alpha_\beta$ to one of the following forms:

\noindent
a) $f^\alpha_\beta
=\delta^\alpha_2\delta_\beta^1$ for $f^2_1 \neq 0$,

\noindent
b) $f^\alpha_\beta = \mu \ve\beta\delta^\alpha_\beta$, $\ve\beta
=(-1)^{\beta+1}$ for $f^2_1= 0$ (here the strikes are ommited).

If we admit the first possibility then come to contradiction with
the assumption that $M_4$ is non-Einstein manifold (see proof in
Appendix~B).

In the second case ($c^2_1=0$) we have
    $$
a^\alpha_\beta = \delta^\alpha_1 \pl\beta \Phi+
\mu \ve\beta \delta^\alpha_\beta + \rho \delta^\alpha_\beta,
\qquad \ve\beta = (-1)^{\beta+1},
\qquad
\mu=\mu (z^2)
    $$
and from Eq.~(\ref{zvez}) it follows
    \bqq{sh52}
R^2_1 = 0, \qquad (\pl 1\Phi +2\mu) R^1_2 + \pl 2 \Phi (R^2_2-R^1_1)=0.
    \eeq
Using the symmetry and reality of $a$, we find
    \bqq{sh53}
g_{1\ol 1}(\mu -\ol\mu)=0,
    \eeq
    $$
g_{2\ol 1}\pl{\ol 2}\Phi - g_{1\ol 2}\pl 2 \Phi = g_{2\ol
2}(\ol\mu-\mu),
    $$
    \bqq{sh54}
g_{1\ol 1}\pl{\ol 2}\Phi - g_{1\ol 2}\pl 1 \Phi = g_{1\ol
2}(\ol\mu+\mu).
    \eeq
If $g_{1\ol 1}\ne 0$ then from (\ref{sh53}) it follows that $\mu=\ol\mu$.
In the case $g_{1\ol 1} =g_{2\ol 2}= 0$ from (\ref{sh54}) and
non-degeneracy of $g$ we find $\pl 1 \Phi =-(\mu+\ol\mu)$. Similarly,
using (\ref{zamkn}) and (\ref{kahl}), we have
     \bqq{sh55}
\pl 1 g_{1\ol 2}=\pl 1 g_{2\ol 1}=\pl 1 g_{2\ol 2}=0.
     \eeq
Therefore, $R^1_1=R^2_2=0$ (see (\ref{ric})) and from (\ref{sh52})
$R^1_2 (\mu - \ol \mu)=0$. From here $R^1_2=0$ for $\mu \neq \ol\mu$,
whence, $R^i_j=0$ and $M_4$ is Einstein manifold that contradict to
our initial assumption. Therefore, $\mu=\ol\mu=$const.

Making the transformations $\rho'=\rho-\mu$, $\Phi'=\Phi+2\mu(z^1
+\ol{z^1})$ and ommiting strikes, we can reduce $a^\alpha_\beta$ to
the form
    \bqq{sh56}
a^\alpha_\beta = \delta_1^\alpha \pl\beta\Phi
+\rho\delta^\alpha_\beta.
    \eeq
The reality and symmetry conditions $g_{2\ol 1}\pl{\ol 2}\Phi -
g_{1\ol 2}\pl 2 \Phi = 0$, $g_{1\ol 1}\pl{\ol 2}\Phi - g_{1\ol 2}\pl 1
\Phi = 0$ of the tensor $a$ are equivalent to the equation
    \bqq{sh57}
\pl 2\Phi =\varphi\pl 1\Phi
    \eeq
where $\varphi=\varphi(z^2,\ol{z^2})$ is a (complex) function
depending on $x^2 = \frac{1}{\sqrt{2}} (z^2+\ol{z^2})$ and
$y^2=\frac{1}{i\sqrt{2}}(z^2-\ol {z^2})$. From here with the use of
Eq.~(\ref{kahl}) we find
    $$
g_{2\ol 1} = \varphi g_{1\ol 1}, \quad g_{1\ol 2}= \ol\varphi g_{1\ol
1},\quad g_{2\ol 2}= \pl{\ol 2}\varphi \pl 1\Phi+ \varphi\ol\varphi
g_{1\ol 1}.
    $$
>From the condition $g_{2\ol 2}= \ol{g_{2\ol 2}}$ it follows $\pl{\ol
2}\varphi\pl 1\Phi=\pl 2\ol\varphi\pl{\ol 1}\Phi$, and because of
Lemma~\ref{sh-lm1} and the formula $\pl 1 \Phi\neq 0$
    \bqq{sh58}
\pl{\ol 2}\varphi= \pl 2\ol\varphi.
    \eeq
This equation can be treated as integrability condition of the system
    \bqq{sysur}
\varphi = \pl 2 F, \qquad \ol\varphi = \pl{\ol 2} F
    \eeq
where $F$ is a real function depending on $z^2$ and $z^{\ol 2}$. If
the equation (\ref{sh58}) holds, then (\ref{sysur}) has a solution
$F$. Substituting it in (\ref{sh57}), we find
    \bqq{sh59}
\pl 2\Phi= \pl 2 F\pl 1 \Phi.
    \eeq
Here $\pl 2\pl {\ol 2} \Phi\neq 0$, because in opposite case $\det \:
(g_{\alpha\ol\beta})= \pl 1 \pl 1 \Phi\pl 1\Phi\pl{\ol 2}\varphi =0$.

The Eq.~(\ref{hpa}) holds identically because of (\ref{sh56}),
(\ref{sh57}) and (\ref{sh59}). Let $\tilde F$ be a real function
functionally independent from $F$. Introducing new variables
$u=F(z^2,\ol{z^2})$ and $v=\tilde F (z^2,\ol{z^2})$, we obtain from
(\ref{sh59})
      $$
\pl u\Phi+\frac{\pl 2 \tilde F}{\pl 2 F}\pl v\Phi=\pl 1\Phi.
      $$
>From here, taking into account reality of the functions $F$, $u$
and $v$ as well as the identity $\dmin\Phi=0$, we find
      $$
(\frac{\pl 2\tilde F}{\pl 2 F}- \frac{\pl{\ol{2}}\tilde F}{\pl{\ol 2}
F})\pl{v}\Phi=0.
      $$
Since $F$ and $\tilde F$ are functionally independent $\pl v\Phi=0$.
Therefore, $\Phi =\Phi(z^1+\ol{z^1} +F)$.

>From these relations the main result now follows
      \begin{The} \label{sh-th2}
Let $f$ be non-affine $H$-projective mapping of non-Einstein
four-dimensional K\"ahler manifold $(M_4,g,J)$ on K\"ahler manifold
$(M'_4,g',J)$. Then in a neighborhood of any point $p\in M_4$ there
exist complex coordinates $(z^{\alpha},z^{\ol\alpha})$,
$\alpha=1,\ldots,n$ in which K\"ahler potential $\Phi$ and
components of the metric $g$ are defined by the formulas
     \bqq{sh60}
\Phi ={\cal W} (z^1+\ol{z^1}+F(z^2,\ol{z^2})), \quad F=\ol F, \: \pl 2
F\neq 0,\quad {\cal W}\neq {\rm const},
     \eeq
     \bqq{kahl-a}
g_{\alpha\ol\beta}=\partial_{\alpha}\partial_{\ol\beta}\Phi.
     \eeq

In corresponding coordinates in a neighborhood of the point $f(p) \in
M'_{4}$ components of the metric $g'$ are defined by Eq.~{\rm
(\ref{sin})} where
     \bqq{sh69}
a_{\alpha\ol\beta}= \ol{a_{\ol\alpha\beta}} = \pl\alpha\Phi\pl
1\pl{\ol\beta}\Phi+ \rho\pl\alpha\pl{\ol\beta} \Phi,\quad
a_{\alpha\beta}=a_{\ol\alpha\ol\beta}=0, \quad \rho\in\rr.
     \eeq
     \end{The}


\section{Generalized equidistant K\"ahler manifolds}


If in (\ref{sh60}) ${\cal W} =\exp (z^1+\ol{z^1}+F(z^2,\ol{z^2}))$ then
(\ref{kahl-a}) defines the metric of an equidistant K\"ahler
manifold \cite{skm}. J.~Mike\v s \cite{mik-ekv} proved that equidistant
K\"ahler manifolds admit non-affine $H$-projective mappings and
Theorem~\ref{sh-th2} confirms this result. The K\"ahler manifolds with
potential (\ref{sh60}) will be called {\it generalized equidistant
K\"ahler manifolds}.

Let us consider a K\"ahler metric $g$ with K\"ahler potential
(\ref{sh60}), the tensor field $a$, defined by Eqs.~(\ref{sh69}) and
gradient covector $\lambda_\alpha = g_{\alpha \ol 1}$,
$\lambda_{\ol\alpha} = g_{\ol\alpha 1}$. As it was shown in the
previous section, Eq.~(\ref{hpa}) holds identically for $g$, $a$ and
$\lambda_\alpha$. Therefore, any generalized equidistant K\"ahler
manifold admits non-affine $H$-projective mapping. Einstein
generalized equidistant K\"ahler manifolds are distinguished by the
condition
        $$
\exp (a\Phi)\pl 1 (\pl 1 \Phi)^2 \pl 2\pl {\ol 2} F = f(z) \ol{f(z)},
\qquad a\in \rr
        $$
where $f(z)$ is an appropriate holomorphic function. Four-dimensional
K\"ahler manifolds of constant holomorphic sectional curvature are
distinguished by
        $$
{\cal W} = \ln (1 + \exp (z^1+\ol{z^1}+ \ln (1+ \epsilon z^2
\ol{z^2}))), \qquad \epsilon = \pm 1.
        $$

Let us now consider Ricci-flat generalized equidistant K\"ahler
manifolds. Ricci-flat K\"ahler manifolds possess the hyper K\"ahler
structure of are and have important applications in theoretical
physics \cite{bess,chsw1,chsw2,azk}.

Rewritting the condition $R_{ij}=0$ with the use of (\ref{ric}), we
obtain
        $$
\det (g_{\mu\ol\nu}) = f(z) \ol{f(z)}
        $$
where $f$ is a holomorphic function. Substituting (\ref{sh60})
in this formula, we find
        $$
{\cal W}' {\cal W}'' \pl{2\ol 2} F=f(z)\ol{f(z)}.
     $$
Eq.~(\ref{sh60}) together with the Note~\ref{note1} at
p.~\pageref{note1} yields $f=$const. Therefore, the condition
$R_{ij}=0$ can be written in the form
     \bqq{eq-W}
{\cal W}' {\cal W}'' ={\rm const}, \qquad \pl{2\ol 2} F={\rm const}.
     \eeq
The last equation gives
     \bqq{z2}
F(z^2,\ol{z^2})=\gamma z^2 \ol{z^2}+\tau (z^2+\ol{z^2})+\rho
     \eeq
where $\gamma,\tau,\rho=$const. Integrating (\ref{eq-W}), we find
     \bqq{eq-W2}
{\cal W}=A(x+B)^{3/2}+C, \qquad
x=z^1+\ol{z^1}+F(z^2,\ol{z^2})
     \eeq
where $A,B$ and $C$ are some constants. After substituting (\ref{z2})
in (\ref{eq-W2}) and changing the variables $z^1\to
z^1+(\tau^2-\rho)/2$, $z^2\to z^2-\tau$, we find the following general
expression for K\"ahler potential of four-dimensional generalized
equidistant manifold
     $$
\Phi=A(z^1+\ol{z^1}+\gamma z^2\ol{z^2})^{3/2}.
     $$
Here the non-essential constant $C$ is omitted. From the last formula
it is easy to find the metric
     $$
ds^2=\frac{3}{4}A (z^1+\ol{z^1}+\gamma z^2\ol{z^2})^{-1/2}
[dz^1 d\ol{z^1} + \gamma z^2 dz^1 d\ol{z^2} +
\gamma\ol{z^2}dz^2 d\ol{z^1} +
2\gamma(z^1+\ol{z^1}+\frac{3}{2} z^2\ol{z^2})dz^2 d\ol{z^2}].
     $$

A real holomorphic vector field $X$ on a K\"ahler manifold $(M,g,J)$
is called {\it infinitesimal $H$-projective transformation} or {\it
$H$-projective motion} if
     $$
L_X \omega =\varphi\:\iota+\iota\:\varphi - \varphi\circ J\:\iota\circ
J - \iota\circ J \: \varphi\circ J.
     $$
Here $\omega$ is the Riemannian connection of $(M,g)$, $\iota$ denotes
the internal product ($\iota:TM \times\land^m T^*M \to
\land^{m-1}T^*M$, $\iota_Y \theta\equiv \theta (Y,\ldots)$, $Y\in TM$,
$\theta \in \land^p T^*M$) and $\varphi$ is an 1-form on $M$. An
$H$-projective motion is called {\it non-affine} if $\varphi \neq 0$.

K.~Yano and K.~Hiramatu \cite{yahi1} proved that compact Einstein
K\"ahler manifolds admit non-affine $H$-projective motions iff they
have constant holomorphic sectional curvature \cite{kobn}. A K\"ahler
manifold admiting non-affine $H$-projective motion also admits
non-affine $H$-projective mapping \cite{sin1}. With the help of this
fact the Theorem~\ref{sh-th2} yields the following
     \begin{The} \label{sh-th3}
If a compact four-dimensional K\"ahler manifold admits infinitesimal
properly $H$-projective transformations, then it is generalized
equidistant K\"ahler manifold.
     \end{The}


\section*{Acknowledgments}


Author is thankful to Prof.~A.~Aminova and Prof.~J.~Mike\v s for
comments, discussions and suggestions. The work was partially
supported by Russian Foundation for Fundamental Investigations, Grant
No~96-0101031.


\setcounter{section}{0}
\renewcommand{\theequation}{\Alph{section}.\arabic{equation}}
\appsection


Here we provide the proof of Lemma~\ref{sh-lm2}.

It follows from (\ref{hpab-sol}) that $b^\alpha_\beta$ depend only on
$z^1$, $z^2$. Holomorphic coordinate transformations
$z'^\alpha=z'^\alpha (z)$ does not change this result and can be used
for vanishing of $b^1_2$.

Let $b^1_2=0$, consider the following three possibilities in
Eq.~(\ref{zvez-bb}).

\noindent
1) Let $b^2_1=0$ and the tensor field $b$ does not vanishes. Then either
$b^1_1=b^2_2=0$ that contradicts with the assumption about not
vanishing of tensor $b^i_j$ or, because $M_4$ is non-Einstein
manifold, $R^1_2=R^2_1=0$ and $R^1_1 \neq R^2_2$. In the last case it
is possible to find such (in general, complex) functions $v_1$ and
$v_2$ that
   \bqq{sh39-bis}
b^\alpha_\beta = v_1 R^\alpha_\beta + v_2 \delta_\beta^\alpha.
   \eeq

\noindent
2) If $b^2_1 \neq 0$ and $b^1_1 \neq 0$, then from (\ref{zvez-bb}) we have
   $$
R^1_2=0,\qquad \frac{R^1_1 -R^2_2}{2 b^1_1}= \frac{R^2_1}{b^2_1},
   $$
hence, $b^1_1-b^2_2=v_1(R^1_1-R^2_2)$, $b^2_1 = v_1 R^2_1=0$, $b^1_2=
v_1 R^1_2$ for some function $v_1$. By putting $b^\alpha_\beta = v_1
R^\alpha_\beta + \tilde b_\beta^\alpha $, we find $\tilde b^1_1 -
\tilde b^2_2 = \tilde b^1_2= \tilde b^2_1 =0 $ or $\tilde
b^\alpha_\beta = v_2\delta^\alpha_\beta$ where $v_2$ is some function
in the region $U\subset M_4$.

\noindent
3) At last, in the case $b^2_1\neq 0$, $b^1_1=b^2_2=0$ we obtain
$R^1_2=R^1_1-R^2_2=0$. Therefore, it is possible to find such
functions $v_1$, $v_2$ that Eq.~(\ref{sh39-bis}) holds. Putting
$v_1=v_2=0$ in Eq.~(\ref{sh39-bis}), we find $b^i_j=0$. We come to the
conclusion that the formula (\ref{sh39-bis}) describes all possible
cases. In the similar way the relations
      \bqq{sh39-1}
b^{\ol\alpha}_{\ol\beta} = \ol{v_1} R^{\ol\alpha}_{\ol\beta} +
\ol{v_2} \delta_\beta^\alpha, \qquad b^{\ol\alpha}_\beta = v_1
R^{\ol\alpha}_\beta + v_2 \delta_\beta^{\ol\alpha} \equiv 0.
      \eeq
can be obtained. From here because of reality and symmetry of the
tensora $a$, $b$ and Ricci tensor it follows that $v_1$ and $v_2$ are
real-valued functions, i.e. (\ref{sh39-bis}), (\ref{sh39-1}) can be
rewritten in the form of the tensor relation
      $$
b^i_j = v_1 R^i_j + v_2 \delta_j^i.
      $$
>From here we find $v_2 = - (v_1 R)/2n$, whence
    \bqq{sh40}
b^\alpha_\beta = v_1 (R^\alpha_\beta - \frac{R}{2n}
\delta_\beta^\alpha),
    \eeq
    $$
b^\alpha_{\beta,j} = {v_1}_{,j}(R^\alpha_\beta -
\frac{R}{2n}\delta_\beta^\alpha) + v_1 (R^\alpha_\beta - \frac{R}{2n}
\delta_\beta^\alpha)_{,j}.
    $$
Because of (\ref{cris}), (\ref{hpa-mod3}) and (\ref{hpab-sol}) we have
$b^\alpha_{\beta,j} = 0$ and, denoting ${\cal A}=\ln v_1$, we find
    $$
{\cal A}_{,j}=-\frac{(R^\alpha_\beta-\delta^\alpha_\beta R/2n)_{,j}}
{R^\alpha_\beta - \delta^\alpha_\beta R/2n}.
    $$
The right part of this relation doesn't depend on the variable
$y^1=\frac{1}{\sqrt{2}}(z^1 -z^{\ol 1})$, hence, its left part also
doesn't depend on $y^1$. Because of reality of ${\cal A}$ we have
    $$
{\cal A} =\tilde f (z^1 +z^{\ol 1},z^2,z^{\ol 2}) +i\tilde \tau\cdot
(z^1-z^{\ol 1}), \qquad \tilde \tau \in\rr,
    $$
    $$
v_1 = \exp {\cal A} = f (z^1 +z^{\ol 1},z^2,z^{\ol 2})\exp (i
\tau\cdot (z^1 - z^{\ol 1})),
   \qquad
\tau \in\rr.
   $$
Then from (\ref{sh40}) we obtain $b^\alpha_\beta = \exp (2i\tau z^1)
\tilde c^\alpha_\beta(z^2)$. Taking into account (\ref{hpab-sol}), we
find
   \bqq{sh41}
f^\alpha_\beta = -i\int b^\alpha_\beta dz^1 = \exp (2i\tau z^1)
c^\alpha_\beta(z^2) +d^\alpha_\beta(z^2).
    \eeq
Since $f^\alpha_\alpha=0$ we have $c^\alpha_\alpha = d^\alpha_\alpha =0$.

>From (\ref{zvez}), (\ref{hpa-sol2}) and (\ref{sh41})
    \bqq{sh43}
\delta^\alpha_1 \pl\mu \Phi R^\mu_\beta -\pl\beta\Phi R^\alpha_1
+d^\alpha_\mu R^\mu_\beta - d^\mu_\beta R^\alpha_\mu = 0,
    \eeq
    \bqq{sh44}
c^\alpha_\mu R^\mu_\beta - c^\mu_\beta R^\alpha_\mu = 0.
    \eeq

Let us consider separately the following two cases: $c^2_1 \neq 0$ and
$c^2_1 =0$. Under the admissible cooordinate transformations
(\ref{free1}) not changing the form of $\Lambda = \pl1 +\pl{\ol 1}$
the components of (1,1)-ten\-sor field $t$ change as following
    $$
t^{1'}_{1'} =t^1_1 + \frac{dl}{dz^2} t^2_1,
    $$
    $$
t^{1'}_{2'} =(\frac{dm}{dz^2})^{-1} (-\frac{dl}{dz^2} t^1_1 -
(\frac{dl}{dz^2})^2 t^2_1 + \frac{dl}{dz^2} t^2_2 +t^1_2),
    $$
    $$
t^{2'}_{1'} = \frac{dm}{dz^2} t^2_1, \qquad t^{2'}_{2'} =
-\frac{dl}{dz^2} t^2_1 + t^2_2.
    $$

>From here it is seen that if $c^2_1 \neq 0$ then the components
$c^\alpha_\beta$ of tensor field $c$ can be reduced to the form
   $$
c^{\alpha}_{\beta} = \delta^{\alpha}_1 \delta^2_{\beta} \phi+
\delta^{\alpha}_2 \delta^1_{\beta}
   $$
where strikes are ommited and $\phi$ is a holomorphic function
depending on $z^2$. Substituting this expression in (\ref{sh44}), we
find $\phi R^{2}_{1} =R^{1}_{2}$, $R^{1}_{1} =R^{2}_{2}$, hence, with
the help of (\ref{sh43}) we have
   $$
(\pl{1} \Phi +2d^{1}_{1}) R^{1}_{2} =0, \qquad (\pl{1} \Phi
+2d^{1}_{1}) R^{2}_{1} =0.
   $$
If $R^1_2 \neq 0$ or $R^2_1 \neq 0$, then $\pl{1} \Phi = -2d^1_1$ is
holomorphic function and $g_{1\ol 1}=g_{1\ol 2}=0$, whence, the metric
is degenerate. Therefore, $R^1_2 = R^2_1 = R^1_1 - R^2_2 = 0$ that
contradicts with our assumption that $M_4$ is non-Einstein manifold.

So we have $c^2_1=0$. In this case with the use of admissible
coordinate transformations the tensor fields $c$ and $d$ can be
reduced to the form
   \bqq{sh46}
c^\alpha_\beta = \ol{c^{\ol\alpha}_{\ol\beta}} = c^1_1
\delta_\beta^\alpha \ve\beta, \qquad \ve\beta =(-1)^{\beta+1},
   \eeq
   $$
d^\alpha_\beta = \ol{d^{\ol\alpha}_{\ol\beta}} = d^1_1
\delta_\beta^\alpha \ve\beta + \gamma (z^2) \delta^\alpha_1
\delta^2_\beta + \zeta \delta^\alpha_2 \delta^1_\beta
   $$
where $\zeta =0,1$ and strikes are omitted. After admissible gauge
transformation $\Phi' = \Phi +\int \gamma(z^2) dz^2 + \ol{\int
\gamma(z^2) d z^2}$, taking into account (\ref{hpa-sol2}),
(\ref{sh41}), we obtain
   \bqq{sh47}
d^\alpha_\beta = \ol{d^{\ol\alpha}_{\ol\beta}} = d^1_1
\delta_\beta^\alpha\ve\beta+\zeta\delta^\alpha_2 \delta^1_\beta.
   \eeq
Substituting (\ref{sh46}) into (\ref{sh44}), we find $c^1_1 R^1_2 =0$,
$c^1_1 R^2_1 =0$, whence, $c^1_1=0$ or $R^1_2 = R_1^2 =0$. In the last
case from (\ref{sh43}) and (\ref{sh47}) we have $\pl 2 \Phi
(R^2_2-R_1^1)=0$. Since $\pl 2\Phi \neq 0$, we find $R^2_2-R_1^1=0$.
We come to Einstein manifold again, therefore, $c_1^1 = 0$ and
$f^\alpha_\beta = d^\alpha_\beta(z^2)$. After substituting this result
into (\ref{hpa-sol2}), we prove the Lemma~\ref{sh-lm2}.

\appsection

Here we prove that the condition $f^\alpha_\beta =\delta^\alpha_2
\delta^1_\beta$ contradicts with the assumption that considered
K\"ahler manifold is non-Einstein.

>From (\ref{zvez}) and (\ref{sh48}) we find
    \bqq{949}
R^1_2=\pl 2 \Phi R^2_1,  \quad R^1_2 \pl 1 \Phi + \pl 2 \Phi (R^2_2 -
R^1_1)=0.
    \eeq
Writing down the symmetry conditions $a_{\alpha\ol\beta} =a^\mu_\alpha
g_{\mu\ol\beta}=a_{\ol\beta\alpha} = a^{\ol\mu}_{\ol\beta}
g_{\alpha\ol\mu}$ of the tensor field $a$ we obtain with the help of
(\ref{sh48}) the next formulas
    $$
g_{2\ol 1} =g_{1\ol 2}, \qquad g_{1\ol 1}\pl 2 \Phi = g_{2\ol
1}\pl{\ol 1}\Phi + g_{2\ol 2},
    $$
    $$
g_{1\ol 2} \pl 1 \Phi + g_{2\ol 2} = g_{1\ol 1} \pl 2 \Phi,
    $$
    $$
g_{1\ol 2}\pl 2 \Phi - g_{2\ol 1}\pl{\ol 2}\Phi \equiv g_{1\ol 2}(\pl
2 \Phi - \pl{\ol 2}\Phi) =0.
    $$
>From the last equation it follows that either $\pl 2 \Phi =\pl{\ol 2}
\Phi$ or $g_{1\ol 2} = g_{2\ol 1} = 0$.

Let us first take $g_{2\ol 1} =g_{1\ol 2}= 0$, then from (\ref{zamkn})
the equality $\pl 1 g_{2\ol 2} = \pl 2 g_{1\ol 1} = 0$ follows, hence
$\pl 1\pl{\ol 2}\det (g_{\alpha\ol\beta})=0$ and, because of
(\ref{ric}) we find $R_{1\ol 2}=R_{2\ol 1}=0$, therefore,
$R_1^2=R^1_2=0$. Taking into account $\pl 2\Phi \neq 0$,
from (\ref{949}) we obtain $R^1_1-R_2^2=0$, which means that the
manifold $M_4$ is Einstein. We came to contradiction with our initial
assumption. Hence, in addition to the formula $\pl 1\Phi = \pl{\ol
1}\Phi$ we have $\pl 2\Phi = \pl{\ol 2}\Phi$. From here using
Eqs.~(\ref{kahl}) -- (\ref{ric}), it is possible to deduce that all
components of the metric tensor, Christoffel symbols and curvature
tensor are real. Then (\ref{zvez}) can be written as
    $$
R_{\alpha\ol\sigma}a^{\ol\sigma}_{\ol\beta} -
R_{\sigma\ol\beta}a^\sigma_\alpha = 0.
    $$
>From here, putting $\alpha,\beta=1,2$ and using the identities $\pl 2
\Phi\neq 0$, $a^\alpha_\beta=a^{\ol\alpha}_{\ol\beta}$ and
$R_{\alpha\ol\beta}=R_{\ol\alpha\beta}$, we find
$R_{\alpha\ol\beta}=0$, hence, $M_4$ is Ricci-flat that contradicts
with the assumption that $M_4$ is non-Einstein manifold. {\it Q.E.D.}


\end{document}